\documentclass[a4,11pt]{amsart}
\usepackage{amssymb,amsmath,amsthm,txfonts}
\usepackage{cite}

\usepackage{color}

\newtheorem{thm}{Theorem}[section]
\newtheorem{lem}[thm]{Lemma}

\newtheorem{prop}[thm]{Proposition}
\theoremstyle{definition}

\theoremstyle{remark}
\newtheorem{rem}[thm]{\textbf{Remark}}
\newtheorem{rems}[thm]{\textbf{Remarks}}

 \makeatletter
    
    \@addtoreset{equation}{section}
  \makeatother

\makeatletter
      \def\@makefnmark{%
         \leavevmode
            \raise.9ex\hbox{\check@mathfonts
                \fontsize\sf@size\z@\normalfont%
                            \@thefnmark}%
       }
      \makeatother

\newcommand{\D}{\textrm{div}}

\newcommand{\dd}{\textrm{d}}

\begin{document}

\title[]{On the large time $L^{\infty}$-estimates of the Stokes semigroup in two-dimensional exterior domains}
\author[]{Ken Abe}
\date{}
\address[K. Abe]{Department of Mathematics, Graduate School of Science, Osaka City University, 3-3-138 Sugimoto, Sumiyoshi-ku Osaka, 558-8585, Japan}
\email{kabe@sci.osaka-cu.ac.jp}

\subjclass[2010]{35Q35, 35K90}
\keywords{Stokes semigroup, $L^{\infty}$-estimates, Stokes paradox}
\date{\today}

\maketitle


\begin{abstract}
We prove that the Stokes semigroup is a bounded analytic semigroup on $L^{\infty}_{\sigma}$ of angle $\pi/2$ for two-dimensional exterior domains. This result is an end point case of the $L^{p}$-boundedness of the semigroup for $p\in (1,\infty)$, established by Borchers and Varnhorn (1993) and an extension of finite time $L^{\infty}$-estimates studied by the author and Giga (2014). The proof is based on the non-existence result of bounded steady flows (the Stokes paradox) and some asymptotic formula for the net force of the Stokes resolvent.
\end{abstract}

\vspace{15pt}

\section{Introduction}

\vspace{10pt}
We consider the Stokes equations:

\begin{equation*}
\begin{aligned}
\partial_t v-\Delta{v}+\nabla{q}= 0,\quad \D\ v&=0  \qquad \textrm{in}\ \Omega\times (0,\infty),  \\
v &=0\qquad \textrm{on}\ \partial\Omega\times (0,\infty), \\
v&=v_0\hspace{18pt} \textrm{on}\ \Omega\times\{t=0\},
\end{aligned}
\tag{1.1}
\end{equation*}\\
for exterior domains $\Omega\subset \mathbb{R}^{n}$, $n\geq 2$. It is well known that the solution operator (called the Stokes semigroup) 

\begin{align*}
S(t):v_0\longmapsto v(\cdot,t),
\end{align*}\\
forms an analytic semigroup on $L^{p}_{\sigma}$ for $p\in (1,\infty)$, of angle $\pi/2$ \cite{Sl76}, \cite{G81}, i.e. $S(t)v_0$ is a holomorphic function in the half plane $\{ \textrm{Re}\ t>0\}$ on $L^{p}_{\sigma}$. Here, $L^{p}_{\sigma}$ denotes the $L^{p}$-closure of $C_{c,\sigma}^{\infty}$, the space of all smooth solenoidal vector fields with compact support in $\Omega$. The Stokes semigroup $S(t)$ is defined by the Dunford integral of the resolvent of the Stokes operator $A=\mathbb{P}\Delta$ for the Helnholtz projection operator $\mathbb{P}: L^{p}\longrightarrow L^{p}_{\sigma}$ \cite{FM77}, \cite{Miyakawa82}, \cite{SiSo}. See, e.g. \cite{Lunardi} for analytic semigroups.

We say that an analytic semigroup on a Banach space is a \textit{bounded} analytic semigroup of angle $\pi/2$ if the semigroup is bounded in the sector $\Sigma_{\theta}=\{t\in \mathbb{C}\backslash \{0\}\ |\ |\arg{t}|<\theta \}$ for each $\theta\in (0,\pi/2)$. See, e.g. \cite[Definition 3.7.3]{ABHN}. The boundedness in the sector implies the bounds on the positive real line

\begin{align*}
||S(t)||\leq C,\quad ||AS(t)||\leq \frac{C}{t},\quad t>0,  \tag{1.2}
\end{align*}\\
where $||\cdot ||$ denotes an operator norm on a Banach space and $A$ is a generator. The estimates (1.2) are important to study large time behavior of solutions to (1.1). In terms of the resolvent, the boundedness of $S(t)$ of angle $\pi/2$ is equivalent to the estimate 

\begin{align*}
||(\lambda-A)^{-1}||\leq \frac{C}{|\lambda|},\quad \lambda\in \Sigma_{\theta+\pi/2}.   \tag{1.3}
\end{align*}\\
When $\Omega$ is bounded, the point $\lambda=0$ belongs to the resolvent set of $A=\mathbb{P}\Delta$ and the Stokes semigroup is a bounded analytic semigroup on $L^{p}_{\sigma}$ of angle $\pi/2$ for $p\in (1,\infty)$. For a half space, the boundedness of the semigroup follows from explicit solution formulas \cite{McCracken}, \cite{Ukai}, \cite{BM88}.

The boundedness of the Stokes semigroup on $L^{p}_{\sigma}$ for $p\in (1,\infty)$ have been established for exterior domains in $\mathbb{R}^{n}$ for $n\geq 2$. For $n\geq 3$, the boundedness of $S(t)$ on $L^{p}_{\sigma}$ is proved in \cite{BS} based on the resolvent estimate

\begin{align*}
|\lambda|||v||_{L^{p}}+|\lambda|^{1/2} ||\nabla v||_{L^{p}}+||\nabla^{2} v||_{L^{p}}\leq C|| f||_{L^{p}},\quad 1<p<\frac{n}{2},   \tag{1.4}
\end{align*}\\
for $v=(\lambda-A)^{-1}f$ and $\lambda\in \Sigma_{\theta+\pi/2}\cup\{0\}$. The estimate (1.4) implies (1.3) for $p\in (1,n/2)$ and the case $p\in [n/2, \infty)$ follows from a duality. Due to the restriction on $p$, the two-dimensional case is more involved. Indeed, the estimate $||\nabla^{2}v||_{L^{p}}\leq C||Av||_{L^{p}}$ for $p\in [n/2,\infty)$ does not hold \cite{BM92}. For $n=2$, the boundedness of the Stokes semigroup on $L^{p}_{\sigma}$ is proved in \cite{BV} based on layer potentials for the Stokes resolvent.

Recently, the case $p=\infty$ has been developed. When $\Omega$ is a half space, $S(t)$ forms a bounded analytic semigroup on $L^{\infty}_{\sigma}$ of angle $\pi/2$ \cite{DHP}, \cite{Sl03}. For a half space and domains with compact boundaries, we define $L^{\infty}_{\sigma}$ by

\begin{align*}
L^{\infty}_{\sigma}(\Omega)=\left\{f\in L^{\infty}(\Omega)\ \Big|\  \D\ f=0\ \textrm{in}\ \Omega,\ f\cdot N=0\ \textrm{on}\ \partial\Omega \right\}.
\end{align*}\\
Here, $N$ denotes the unit outward normal vector field on $\partial\Omega$. Since $S(t)$ is bounded on $L^{\infty}_{\sigma}$, the associated generator $A=A_{\infty}$ is also defined for $p=\infty$. For bounded domains \cite{AG1} and exterior domains \cite{AG2}, analyticity of the semigroup on $L^{\infty}_{\sigma}$ follows from the a priori estimate

\begin{align*}
||v||_{L^{\infty}}+t^{1/2}||\nabla v||_{L^{\infty}}+t||\nabla^{2}v||_{L^{\infty}}+t||\partial_t v||_{L^{\infty}}+t||\nabla q||_{L^{\infty}}\leq C||v_0||_{L^{\infty}},  \tag{1.5}
\end{align*}\\
for $v=S(t)v_0$ and $t\leq T$. The estimate (1.5) is proved by a blow-up argument and implies that $S(t)$ is analytic on $L^{\infty}_{\sigma}$. Moreover, by the resolvent estimates on $L^{\infty}_{\sigma}$ \cite{AGH}, $S(t)$ is analytic on $L^{\infty}_{\sigma}$ of angle $\pi/2$. When $\Omega$ is bounded, $S(t)$ is a bounded analytic semigroup on $L^{\infty}_{\sigma}$ of angle $\pi/2$.

In this paper, we consider the boundedness of the Stokes semigroup on $L^{\infty}_{\sigma}$ for exterior domains in $\mathbb{R}^{n}$ for $n\geq 2$. For the Laplace operator or uniformly elliptic operators, a standard approach to prove large time $L^{\infty}$-estimates of a semigroup is to use a Gaussian upper bound for a complex time heat kernel. See \cite[Chapter 3]{Davies}. However, a kernel of the Stokes semigroup does not satisfy a Gaussian bound since $S(t)$ is unbounded on $L^{1}$. See \cite{DHP}, \cite{Saal07} for a half space. Even for exterior domains, $S(t)$ is not bounded on $L^{1}$ unless the net force vanishes \cite{Kozono98}, \cite{HeMiyakawa06}. It seems no general method to estimate the $L^{\infty}$-norm of a semigroup for all time without a Gaussian bound.

There is a work by Maremonti \cite{Mar14} who proved the estimate

\begin{align*}
||S(t)v_0||_{L^{\infty}}\leq C||v_0||_{L^{\infty}},\quad t>0,  \tag{1.6}
\end{align*} \\
for exterior domains and $n\geq 3$ based on the finite time estimate in  \cite{AG1}. Subsequently, Hieber and Maremonti \cite{MH16} proved the estimate $t||AS(t)v_0||_{L^{\infty}}\leq C||v_0||_{L^{\infty}}$ for $t>0$ and the results are extended in \cite{BH} for complex time $t\in \Sigma_{\theta}$ and $\theta\in (0,\pi/2)$ based on the approach in \cite{Mar14}. The method in \cite{Mar14} seems a perturbation from the heat equation in $\mathbb{R}^{n}$ and excludes the case $n=2$. 

In the previous work \cite{A8}, the author studied large time $L^{\infty}$-estimates of the Stokes semigroup for $n\geq 2$ based on a Liouville theorem for the Stokes equations introduced by Jia, Seregin and \v{S}ver\'{a}k \cite{JSS}. Liouville theorems are important to study regularity of solutions. See \cite{KNSS}, \cite{SS} for Liouville theorems of the Navier-Stokes equations. They may be also related with large time behavior. Following \cite{JSS}, we say that $v\in L^{1}_{\textrm{loc}}(\overline{\Omega}\times (-\infty,0])$ is an ancient solution to the Stokes equations (1.1) if $\D\ v=0$ in $\Omega\times (-\infty,0)$, $v\cdot N=0$ on $\partial\Omega\times (-\infty,0)$ and 

\begin{align*}
\int_{-\infty}^{0}\int_{\Omega}v\cdot (\partial_t\varphi+\Delta\varphi)\dd x\dd t=0,   
\end{align*}\\
for all $\varphi\in C^{2,1}_{c}(\overline{\Omega}\times (-\infty,0])$ satisfying $\D\ \varphi=0$ in $\Omega\times (-\infty,0)$ and $\varphi=0$ on $\partial\Omega\times (-\infty,0)\cup\Omega\times \{t=0\}$. The conditions $\D\ v=0$ and $v\cdot N=0$ are understood in the sense that 

\begin{align*}
\int_{\Omega}v\cdot \nabla \Phi\dd x=0,\quad \textrm{a.e.}\ t\in (-\infty,0),
\end{align*}\\
for all $\Phi\in C_{c}^{1}(\overline{\Omega})$. Liouville theorems for the Stokes equations has been established in \cite{JSS} for $\mathbb{R}^{n}$, $\mathbb{R}^{n}_{+}$ and bounded domains. Among others, it is proved in \cite{JSS} for exterior domains in $\mathbb{R}^{n}$ for $n\geq 3$ that bounded ancient solutions $v\in L^{\infty}(\Omega\times (-\infty,0))$ must satisfy 

\begin{align*}
v(x,t)-v_{\infty}(t)=O(|x|^{-n+2})\quad \textrm{as}\ |x|\to\infty,
\end{align*}\\
for some constant $v_{\infty}(t)$. Since bounded steady flows exist for $n\geq 3$  \cite{BM92}, bounded ancient solutions are non-trivial. If in addition some spatial decay condition is assumed, we can exclude such solutions.

\vspace{10pt}

\begin{thm}[Liouville theorem on $L^{p}$ \cite{A8}]
Let $\Omega$ be an exterior domain with $C^{3}$-boundary in $\mathbb{R}^{n}$, $n\geq 2$. Let $v$ be an ancient solution to the Stokes equations (1.1). Assume that 

\begin{align*}
v\in L^{\infty}(-\infty,0; L^{p})\quad \textrm{for}\ p\in (1,\infty).  
\end{align*}\\
Then, $v\equiv 0$.
\end{thm}

\vspace{15pt}

Theorem 1.1 is used to prove the large time $L^{\infty}$-estimate (1.6). By the representation formula for $v=S(t)v_0$ \cite{Mizumachi84}, we have

\begin{align*}
v(x,t)=\int_{\Omega}\Gamma (x-y,t)v_0(y)\dd y+\int_{0}^{t}\int_{\partial\Omega}V(x-y,t-s)(TN)(y,s)\dd H(y)\dd s.  \tag{1.7}
\end{align*}\\
Here, $T=\nabla v+{}^{t}\nabla v-qI$ is the stress tensor with the identity matrix $I$ and $V=(V_{ij})$ is the Oseen tensor

\begin{align*}
V_{ij}(x,t)=\delta_{ij}\Gamma(x,t)+\partial_i\partial_j\int_{\mathbb{R}^{n}}E(x-y)\Gamma(y,t)\dd y,   
\end{align*}\\
defined by the heat kernel $\Gamma(x,t)=(4\pi t)^{-n/2}e^{-|x|^{2}/4t}$ and the fundamental solutions of the Laplace equation $E$, i.e.

\begin{align*}
E(x)=
\begin{cases}
&\displaystyle\frac{1}{n(n-2)\alpha(n)}\frac{1}{|x|^{n-2}},\quad n\geq 3, \\
&-\displaystyle\frac{1}{2\pi}\log{|x|},\hspace{47pt} n=2, 
\end{cases}
\end{align*}\\
where $\alpha(n)$ denotes the volume of the unit ball in $\mathbb{R}^{n}$.

For $n\geq 3$, the formula (1.7) describes the asymptotic behavior of bounded Stokes flows as $|x|\to\infty$ and $t\to\infty$. Since the Oseen tensor satisfies 

\begin{align*}
|V(x,t)|\leq \frac{C}{(|x|+t^{1/2})^{n}},\quad x\in \mathbb{R}^{n},\ t>0,  
\end{align*}\\
we have

\begin{align*}
\left|v(x,t)-\int_{\Omega}\Gamma(x-y,t)v_0(y)\dd y\right|
\leq \frac{C}{|x|^{n-2}}\sup_{0<s\leq t}||T||_{L^{\infty}(\partial\Omega)}(s),\quad  |x|\geq R,\ t>0,\tag{1.8}
\end{align*}\\
for some constant $R>0$. The right-hand side is decaying as $|x|\to\infty$ uniformly for all $t>0$. The large time estimate (1.6) for $n\geq 3$ is deduced in \cite{A8} by using the asymptotic formula (1.8) and the Liouville theorem (Theorem 1.1) by a contradiction argument. Indeed, if (1.6) were false, a sequence of solutions generates a non-trivial ancient solution satisfying $|v(x,t)|\leq C|x|^{-n+2}$ for $|x|\geq R$, $t\in (-\infty,0]$ and the Liouville theorem yields a contradiction. The boundedness of $S(t)v_0$ in the sector $\Sigma_{\theta}$ follows the same argument on the half line $\{\arg\ {t}=\theta\}$.

For $n=2$, there is a restriction on the net force since the right-hand side of (1.8) might diverge. Indeed, we have 

\begin{equation*}
\begin{aligned}
&\left|v(x,t)-\int_{\Omega}\Gamma(x-y,t)v_0(y)\dd y-\int_{0}^{t}V(x,t-s)F(s)\dd s\right| \\
&\leq\frac{C}{|x|} \sup_{0<s\leq t}||T||_{L^{\infty}(\partial\Omega)}(s),\quad |x|\geq R,\ t>0,
\end{aligned}
 \tag{1.9}
\end{equation*}
with \textit{the net force}

\begin{align*}
F(s)=\int_{\partial\Omega}TN(y,s)\dd H(y).
\end{align*}\\
Since $|\int_{0}^{t}V(x,s)\dd s|\lesssim \log{(1+t/|x|^{2})}$, the decay as $|x| \to\infty$ of the third term in (1.9) is not uniform for $t>0$ in contrast to (1.8) for $n\geq 3$. If the net force vanishes, the situation is the same as $n=3$ and we are able to prove (1.6) for $t\in \Sigma_{\theta}$. For example, when $\Omega^{c}$ is a disk and initial data has some discrete symmetry (called $C_m$-covariance), the net force vanishes \cite{HeMiyakawa06}, i.e. $F(s)\equiv 0$. The following result includes the case $n=2$ which seems first appeared in \cite{A8}.

\vspace{15pt}

\begin{thm}[Boundedness on $L^{\infty}$ for $n\geq 3$ and $n=2$ with zero net force \cite{A8}]

\noindent
(i) For $n\geq 3$, the Stokes semigroup is a bounded analytic semigroup on $L^{\infty}_{\sigma}$ of angle $\pi/2$. 

\noindent
(ii) For n=2, the estimate (1.6) holds for $t\in \Sigma_{\theta}$ and $v_0\in L^{\infty}_{\sigma}$ for which the net force vanishes (e.g. $C_m$-covariant vector fields when $\Omega^{c}$ is a disk.)
\end{thm}

\vspace{15pt}

In this paper, we prove that the assertion (ii) of Theorem 1.2 holds for any bounded initial data $v_0\in L^{\infty}_{\sigma}$. Perhaps the most important vector fields with non-vanishing net force are asymptotically constant solutions of the steady Navier-Stokes flows as $|x|\to\infty$ such as D-solutions or PR-solutions. See \cite{Gal04}. They are bounded and with finite Dirichlet integral. The situation is subtle even for bounded initial data with finite Dirichlet integral for which the fractional power estimate

\begin{align*}
||\nabla v||_{L^{2}}=||(-A)^{1/2}v||_{L^{2}},
\end{align*}\\
is available. This estimate holds only for $n=2$, i.e.  the estimate $||\nabla v||_{L^{p}}\leq C ||(-A)^{1/2}v||_{L^{p}}$ for $p\in [n,\infty)$ and $n\geq 3$ does not hold \cite{BM92}. The fractional power estimate implies a uniform bound in the homogeneous $L^{2}$-Sobolev space $\dot{H}^{1}$ and $S(t)v_0$ is merely bounded in $\textrm{BMO}$ even if $v_0$ is with finite Dirichlet integral, i.e.

\begin{equation*}
[S(t)v_0]_{\textrm{BMO}}\leq C||v_0||_{L^{\infty}\cap \dot{H}^{1}},\quad t> 0.
\end{equation*}\\
To prove the large time $L^{\infty}$-estimate (1.6) for $n=2$ and any bounded initial data $v_0\in L^{\infty}_{\sigma}$, we analyze the corresponding Stokes resolvent problem:

\begin{equation*}
\begin{aligned}
\lambda v-\Delta v+\nabla q=f,\quad \D\ v&=0\quad \textrm{in}\ \Omega, \\
v&=0\quad \textrm{on}\ \partial\Omega.
\end{aligned}
\tag{1.10}
\end{equation*}\\
Existence and uniqueness of the problem (1.10) for $f\in L^{\infty}_{\sigma}$ have been studied in \cite{AGH}. In particular, the solution operator 

\begin{equation*}
R(\lambda): f\longmapsto v(\cdot,\lambda),
\end{equation*}\\
is a bounded operator on $L^{\infty}_{\sigma}$ for $\lambda\in \Sigma_{\theta+\pi/2}$ and for each $\delta >0$, the estimate $||R(\lambda)||\leq C_{\delta}|\lambda|^{-1}$ holds for $|\lambda|\geq \delta$ with the operator norm $||\cdot||$ on $L^{\infty}_{\sigma}$. The operator $R(\lambda)$ is resolvent of some closed operator $A=A_{\infty}$ on $L^{\infty}_{\sigma}$, i.e. $R(\lambda)=(\lambda-A)^{-1}$. The behavior of $R(\lambda)$ as $\lambda\to0$ corresponds to the behavior of $S(t)$ as $t\to\infty$. Instead of proving the boundedness of $S(t)$ in $\Sigma_{\theta}$, we shall prove the equivalent estimate (1.3) with the operator norm on $L^{\infty}_{\sigma}$. The main result of this paper is the following:

\vspace{10pt}

\begin{thm}[Boundedness on $L^{\infty}$ for $n=2$]
Let $\Omega$ be an exterior domain with $C^{3}$-boundary in $\mathbb{R}^{2}$.

\noindent
(i) For $\theta\in (0,\pi/2)$, there exists a constant $C$ such that 

\begin{align*}
||R(\lambda)f||_{L^{\infty}} \leq \frac{C}{|\lambda|}||f||_{L^{\infty}}, 
\quad \lambda\in \Sigma_{\theta+\pi/2},\ f\in L^{\infty}_{\sigma}.
\tag{1.11}
\end{align*}\\
(ii) The Stokes semigroup is a bounded analytic semigroup on $L^{\infty}_{\sigma}$ of angle $\pi/2$. 
\end{thm}

\vspace{15pt}

There is a difference on the large time behavior for $n=2$ and $n\geq 3$. By Theorems 1.2 and 1.3, we obtain

\begin{align*}
||S(t)v_0||_{L^{\infty}}+t||AS(t)v_0||_{L^{\infty}}\leq C||v_0||_{L^{\infty}},\quad t>0,\ v_0\in L^{\infty}_{\sigma},     \tag{1.12}
\end{align*}\\
for exterior domains in $\mathbb{R}^{n}$ for $n\geq 2$. The estimate (1.12) implies that $S(t)v_0$ is uniformly bounded and approaches a steady flow as $t\to\infty$. For $n=2$, any bounded solutions of 

\begin{equation*}
\begin{aligned}
-\Delta v+\nabla q=0,\quad \D\ v&=0\quad \textrm{in}\ \Omega,\\
v&=0\quad \textrm{on}\ \partial\Omega,
\end{aligned}
\tag{1.13}
\end{equation*}\\
must be trivial (the Stokes paradox) \cite{ChangFinn} and therefore $S(t)v_0$ converges to zero locally uniformly in $\overline{\Omega}$ as $t\to\infty$. On the other hand, for $n\geq 3$, bounded steady flows of (1.13) exist and must be asymptotically constant as $|x|\to\infty$. Hence the solution $S(t)v_0$ converges to such a stationary solution as $t\to\infty$. 

If initial data $v_0$ is decaying as $|x|\to\infty$, $S(t)v_0$ vanishes as $t\to\infty$ for all dimensions $n\geq 2$, i.e. for $v_0\in C_{0,\sigma}$, $S(t)v_0$ uniformly converges to zero in $\overline{\Omega}$ as $t\to\infty$. Here, $C_{0,\sigma}$ is the $L^{\infty}$-closure of $C_{c,\sigma}^{\infty}$, characterized by

\begin{align*}
C_{0,\sigma}(\Omega)
=\left\{f \in C(\overline{\Omega})\ \middle|\ \D\ f=0\ \textrm{in}\ \Omega, \ f=0\ \textrm{on}\ \partial\Omega,\ \lim_{|x|\to\infty}f(x)=0\   \right\}.
\end{align*}\\
See \cite{AG2}. Since $S(t)v_0$ vanishes as $t\to\infty$ for $v_0\in C^{\infty}_{c,\sigma}$, this property follows from the density in $C_{0,\sigma}$. 

There is some issue on the large time behavior of Navier-Stokes flows. By a perturbation argument from the Stokes flow, we are able to construct a unique global-in-time solution of the two-dimensional Navier-Stokes equations for bounded initial data with finite Dirichlet integral \cite{A7} satisfying the integral form

\begin{align*}
u(t)=S(t)u_0-\int_{0}^{t}S(t-s)\mathbb{P}u\cdot \nabla u(s)\dd s.  \tag{1.14}
\end{align*}\\
This solution is asymptotically constant if $u_0$ is, cf. \cite{MaremontiShimizu}. The large time behavior of this solution is an interesting question since the space $L^{\infty}\cap \dot{H}^{1}$ includes steady Navier-Stokes flows. See \cite{Maekawa19} for stability of PR-solutions. It is a question whether solutions of (1.14) remain bounded for all time. The estimate (1.12) implies that the Stokes flow remains bounded for all time and converges to zero locally uniformly in $\overline{\Omega}$ as $t\to\infty$ for any bounded initial data.

The question is non-trivial even for the Cauchy problem for which solutions remain bounded in $\dot{H}^{1}$ by an a priori estimate of vorticity. This solution is merely bounded in BMO. But a uniform $L^{\infty}$-bound seems unknown. The problem have been studied for merely bounded initial data $u_0\in L^{\infty}_{\sigma}$ and a polynomial growth bound on the $L^{\infty}$-norm is derived in \cite{Zelik13}. It is known that global-in-time solutions satisfy the upper bound $||u||_{L^{\infty}}=O(t)$ as $t\to\infty$ \cite{GallayL}. See also \cite{GallaySl}.

\vspace{15pt}

We sketch the proof of Theorem 1.3. Our proof is based on the representation formula for the Stokes resolvent $v=R(\lambda)f$:

\begin{align*}
v(x)=\int_{\Omega}E^{\lambda}(x-y)f(y)\dd y
+\int_{\partial\Omega}V^{\lambda}(x-y)TN(y)\dd H(y), \tag{1.15}
\end{align*}\\
for $T=\nabla v+{}^{t}\nabla v-qI$. Here, 

\begin{align*}
E^{\lambda}(x)=\frac{1}{2\pi}K_{0}(\sqrt{\lambda}|x|)  \tag{1.16}
\end{align*}\\
is the kernel of the resolvent $(\lambda-\Delta)^{-1}$ and $K_m(\kappa)$ is the modified Bessel function of the second kind of order $m$. For $\lambda\in \Sigma_{\theta+\pi/2}$, $\sqrt{\lambda}$ denotes the square-root of $\lambda$ with positive real part, i.e. $\textrm{Re}\ \sqrt{\lambda}>0$. The tensor $V^{\lambda}=(V^{\lambda}_{ij})$ is the kernel of $\lambda (\lambda-\Delta)^{-1}\mathbb{P}$ for the Helmholtz projection operator $\mathbb{P}=I+\nabla(-\Delta)^{-1}\D$. This tensor has the explicit form \cite[p.281]{BV},  

\begin{align*}
V_{ij}^{\lambda}(x)=\frac{1}{2\pi}\left(\delta_{ij} e_{1}\left(\sqrt{\lambda}|x| \right)+\frac{x_ix_j}{|x|^{2}}e_2\left(\sqrt{\lambda}|x|\right)   \right),  \tag{1.17}
\end{align*}\\
where

\begin{align*}
e_1(\kappa)&=K_0(\kappa)+\kappa^{-1}K_1(\kappa)-\kappa^{-2},\\
e_2(\kappa)&=-K_0(\kappa)-2\kappa^{-1}K_1(\kappa)+2\kappa^{-2},\quad \kappa>0.
\end{align*}\\
The function $e_1(\kappa)$ has a logarithmic singularity as $\kappa\to0$ and decaying as $\kappa\to\infty$. The function $e_2(\kappa)$ is bounded for $\kappa>0$, i.e. 

\begin{equation*}
\begin{aligned}
\left|e_1(\kappa)+\frac{1}{2}\log{\kappa}\right|+|e_2(\kappa)|
&\leq C,\hspace{40pt} 0<\kappa\leq d, \\ 
|e_1(\kappa)|+|e_2(\kappa)|
&\leq C \kappa^{-2},\hspace{42pt} \kappa\geq  d, 
\end{aligned}
\tag{1.18}
\end{equation*}\\
for any $d>0$ with some constant $C$. Hence  

\begin{align*}
V^{\lambda}(x)=-\frac{1}{4\pi}\left(\log{\sqrt{\lambda}}+\log{|x|}   \right)I+\tilde{V}^{\lambda}(x),   \tag{1.19}
\end{align*}\\
with a bounded function $\tilde{V}^{\lambda}$ for $|\lambda|^{1/2}|x|\leq d$. For $|\lambda|^{1/2}|x|\geq d$, $V^{\lambda}$ is bounded.

We shall suppose that $\lambda v$ is uniformly bounded on $L^{\infty}$ and observe the asymptotic behavior of $|\lambda|\ ||v||_{L^{\infty}}$ as $\lambda \to0$. We take a point $x_{\lambda}\in \Omega$ such that 

\begin{align*}
||v||_{L^{\infty}}\approx |v(x_{\lambda})|.
\end{align*}\\
The behavior of $\lambda v$ as $\lambda\to0$ is related with the behavior of $f$ as $|x|\to\infty$. For simplicity of the explanation, we shall consider positive $\lambda>0$ and asymptotically constant vector fields $f\to f_{\infty}$ as $|x|\to\infty$ for which $\lambda (\lambda-\Delta)^{-1}f\to f_{\infty}$ as $\lambda\to 0$.

We first observe that $\lambda v$ converges to zero locally uniformly in $\overline{\Omega}$ as $\lambda\to0$. Indeed, since $u=\lambda v$ is uniformly bounded on $L^{\infty}$ and satisfies 

\begin{equation*}
\begin{aligned}
\lambda u -\Delta u+\nabla p=\lambda f,\quad \D\ u&=0\quad \textrm{in}\ \Omega,\\
u&=0\quad \textrm{on}\ \partial\Omega,
\end{aligned}
\tag{1.20}
\end{equation*}\\
for $p=\lambda q$, by elliptic regularity, $u$ converges to a limit locally uniformly in $\overline{\Omega}$ together with $\nabla u$ and $p$. This pressure $p$ is unique up to constant. Since any bounded solutions of (1.13) must be trivial by the Stokes paradox, it turns out that $u$, $\nabla u$ and $p$ converge to zero locally uniformly in $\overline{\Omega}$. This in particular implies that the stress tensor $T=\nabla u+{}^{t}\nabla u-pI$ vanishes on $\partial\Omega$ as $\lambda\to0$.

The behavior of $|\lambda|\ ||v||_{L^{\infty}}=|u(x_{\lambda})|$ depends on that of the points $\{x_{\lambda}\}$. If the points $\{x_{\lambda}\}$ remain bounded, $u(x_{\lambda})$ converges to zero as $\lambda \to0$, i.e. $\lim_{\lambda\to 0}|u(x_{\lambda})|=0$. If the points $\{x_{\lambda}\}$ diverge, according to the logarithmic singularity of $e_1(\kappa)$ as $\kappa\to 0$, we consider two cases whether $\liminf_{\lambda\to 0}|\lambda|^{1/2}|x_\lambda|>0$ or $\liminf_{\lambda\to 0}|\lambda|^{1/2}|x_\lambda|=0$. If $\liminf_{\lambda\to 0}|\lambda|^{1/2}|x_\lambda|>0$, the kernel $V^{\lambda}(x_{\lambda})$ remains bounded by (1.18). Substituting $x=x_{\lambda}$ into 

\begin{equation*}
\begin{aligned}
u(x)=\lambda (\lambda-\Delta)^{-1}f
+\int_{\partial\Omega}V^{\lambda}(x-y)TN(y)\dd H(y),
\end{aligned}
\tag{1.21}
\end{equation*}\\
and sending $\lambda\to 0$ implies $\limsup_{\lambda \to 0}|u(x_{\lambda})|\leq ||f||_{L^{\infty}}$.

If $\liminf_{\lambda\to 0}|\lambda|^{1/2}|x_\lambda|=0$, the kernel $V^{\lambda}(x_{\lambda})$ can be singular as $\lambda\to0$. By (1.19),  

\begin{equation*}
\begin{aligned}
u(x)&=\lambda (\lambda-\Delta)^{-1}f
-\frac{1}{4\pi}\log\sqrt{\lambda}\int_{\partial\Omega}TN(y)\dd H(y) \\
&-\frac{1}{4\pi}\int_{\partial\Omega}\log |x-y| TN(y)\dd H(y) 
+\int_{\partial\Omega}\tilde{V}^{\lambda}(x-y)TN(y)\dd H(y).  
\end{aligned}
\tag{1.22}
\end{equation*}\\
For fixed $x\in \Omega$, sending $\lambda\to0$ implies \textit{the asymptotic formula for the net force}:

\begin{align*}
0=f_{\infty}-\frac{1}{4\pi}\lim_{\lambda \to 0}\log{\sqrt{\lambda}}\int_{\partial\Omega}TN(y)\dd H(y).   \tag{1.23}
\end{align*}\\
The formula (1.23) has been derived for the Oseen approximation by Finn and Smith \cite{FinnSmith67a}. It implies that the net force is asymptotically pure drag, i.e. the direction of the net force is asymptotically same as the uniform flow $f_{\infty}$ as $\lambda\to0$. By choosing a subsequence, we may assume that $|\lambda|^{1/2}|x_{\lambda}|\to0$. We substitute $x=x_{\lambda}$ into (1.22) and send $\lambda\to 0$. Since $|x_\lambda|\leq |\lambda|^{-1/2}$ for small $\lambda>0$, we have

\begin{align*}
\frac{1}{4\pi}\left|\int_{\partial\Omega}\log |x_{\lambda}-y| TN(y)\dd H(y)\right|
\leq -\frac{1}{4\pi}\log{|\lambda|^{1/2}}\left|\int_{\partial\Omega} TN(y)\dd H(y)\right|+o(1)\quad \textrm{as}\ \lambda\to 0.
\end{align*}\\
By (1.23), $\limsup_{\lambda \to 0}|u(x_{\lambda})|\leq ||f||_{L^{\infty}}$. Hence in all cases, the sup-norm of $\lambda v=u$ is controlled by that of $f$. 

Based on this observation, we apply a contradiction argument to obtain the desired estimate (1.11). We suppose that (1.11) were false and obtain sequences $\{f_m\}$ and $\{\lambda_m\}\subset \Sigma_{\theta+\pi/2}$ such that 

\begin{align*}
&\sup_{\lambda\in \Sigma_{\theta+\pi/2}}|\lambda|\ ||R(\lambda)f_m||_{L^{\infty}}=1,\quad ||f_m||_{L^{\infty}}< \frac{1}{m},\\
&|\lambda_m|\ ||R(\lambda_m)f_m||_{L^{\infty}}\geq \frac{1}{2},\quad \lambda_m\to0.
\end{align*}\\
We set $u_m=\lambda_m R(\lambda_m)f_m$ and take a point $x_m\in \Omega$ such that $|u_m(x_m)|\geq 1/4$. Since $u_m$ satisfies the Stokes resolvent equations (1.20) for $\lambda_m$ with the associated pressure $p_m$, $u_m$ converges to zero locally uniformly in $\overline{\Omega}$ together with $\nabla u_m$ and $p_m$. Then, there are two cases whether $\liminf_{m\to \infty}|\lambda_m|^{1/2}|x_m|>0$ or $\liminf_{m\to \infty}|\lambda_m|^{1/2}|x_m|=0$. Since $||f_m||_{L^{\infty}}\to 0$, in all cases we will see that $1/4\leq |u_m(x_m)|\to 0$ as $m\to\infty$. This is a contradiction.  

\vspace{15pt}
This paper is organized as follows. In Section 2, we prove the representation formula (1.15) for solutions of (1.10) for bounded data $f\in  L^{\infty}_{\sigma}$ and non-existence of bounded solutions of (1.13). In Section 3, we prove Theorem 1.3. After the proof of Theorem 1.3, we note large time behavior of $S(t)v_0$ for $v_0\in L^{\infty}_{\sigma}$.

\vspace{20pt}

\section{Stokes resolvent on $L^{\infty}_{\sigma}$}

\vspace{15pt}

We recall some existence and uniqueness result for the Stokes resolvent equations (1.10) for bounded data $f\in L^{\infty}_{\sigma}$. To state a result, let $L^{p}_{\textrm{ul}}(\overline{\Omega})$ denote the uniformly local $L^{p}$-space in $\overline{\Omega}$ for $p\in (1,\infty)$ and $W^{2,p}_{\textrm{ul}}(\overline{\Omega})$ denote the space of all uniformly local $L^{p}$-functions up to second orders. Let $L^{\infty}_{d}(\Omega)$ denote the space of all functions $f\in L^{1}_{\textrm{loc}}(\overline{\Omega})$ such that $df\in L^{\infty}(\Omega)$ with the distance function $d(x)=\inf_{y\in \partial\Omega}|x-y|$.

\vspace{15pt}

\begin{lem}[Resolvent estimates for large $\lambda$]
(i) For  $p>2$, $\delta>0$ and $\theta\in (0,\pi/2)$, there exists a constant $C$ such that for $f\in L^{\infty}_{\sigma}$ and $\lambda\in \Sigma_{\theta+\pi/2}$ satisfying $|\lambda|\geq \delta$, there exists a unique solution $(v,\nabla q)\in W^{2,p}_{\textrm{ul}}(\overline{\Omega})\times (L^{p}_{\textrm{ul}}(\overline{\Omega})\cap L^{\infty}_{d}(\Omega))$ of (1.10) satisfying 

\begin{align*}
|\lambda| ||v||_{L^{\infty}}+|\lambda|^{1/2}||\nabla v||_{L^{\infty}}+|\lambda|^{1/p}\sup_{x\in \Omega} \left\{ ||\nabla^{2}v||_{L^{p}(\Omega_{x,|\lambda|^{-1/2}}) }+||\nabla q||_{L^{p}(\Omega_{x,|\lambda|^{-1/2}}) } \right\} \leq C||f||_{L^{\infty}},    \tag{2.1}
\end{align*}\\
for $\Omega_{x,r}=\Omega\cap B(x,r)$, where $B(x,r)$ denotes an open ball centered at $x$ with radius $r$.

\noindent
(ii) The solution operator $R(\lambda): f\longmapsto v$ is a bounded operator on $L^{\infty}_{\sigma}$ and satisfies

\begin{align*}
||R(\lambda)||\leq \frac{C}{|\lambda|},\quad \lambda\in \Sigma_{\theta+\pi/2},\    |\lambda|\geq \delta,      \tag{2.2}
\end{align*}\\
with the constant $C$ depending on $\delta$, where $||\cdot ||$ denotes the operator norm on $L^{\infty}_{\sigma}$.
\end{lem}

\vspace{5pt}

\begin{proof}
See \cite[Theorems 1.1 and 1.3]{AGH}.
\end{proof}

\vspace{15pt}
The a priori estimate (2.1) is obtained by applying the localization technique of Masuda \cite{Masuda72} and Stewart \cite{Ste74} by using the $L^{\infty}$-estimate of the pressure. See (2.5) below. The uniqueness follows the same argument. The existence is based on the following approximation lemma for $f\in L^{\infty}_{\sigma}$.

\vspace{15pt}

\begin{lem}[Approximation]
(i) There exists a constant $C$ such that for $f\in L^{\infty}_{\sigma}$ there exists a sequence $\{f_m\}\subset C_{c,\sigma}^{\infty}$ such that 

\begin{equation*}
\begin{aligned}
&||f_m||_{L^{\infty}}\leq C||f||_{L^{\infty}}\\
&f_m\to f\quad \textrm{a.e. in }\ \Omega\quad \textrm{as}\ m\to\infty.
\end{aligned}
\tag{2.3}
\end{equation*}\\
(ii) The resolvent $R(\lambda)f_m$ converges to $R(\lambda )f$ locally uniformly in $\overline{\Omega}$ as $m\to\infty$ for each $\lambda\in \Sigma_{\theta+\pi/2}$.
\end{lem}

\vspace{5pt}

\begin{proof}
The assertion (i) is proved in \cite[Lemma 5.1]{AG2} by using the Bogovski\u\i\ operator. Since $R(\lambda)f_m$ is resolvent of the Stokes operator on $L^{p}_{\sigma}$, i.e. $R(\lambda)f_m=(\lambda-A)^{-1}f_m$ for $A=\mathbb{P}\Delta$, the assertion (ii) follows by applying the a priori estimate (2.1) and uniqueness of (1.10) \cite{AGH}.
\end{proof}

\vspace{15pt}

\begin{rems}
(i) The associated pressure $q$ of the problem (1.10) is a solution of the Neumann problem

\begin{equation*}
\begin{aligned}
-\Delta q&=0\quad \textrm{in}\ \Omega,\\
N\cdot \nabla q&=-N\cdot \nabla^{\perp} \omega\quad \textrm
{on}\ \partial\Omega.
\end{aligned}
\tag{2.4}
\end{equation*}\\
for $\omega=\partial_1v^{2}-\partial_2 v^{1}$, $v={}^{t}(v^{1}, v^{2})$ and $\nabla^{\perp}={}^{t}(\partial_2,-\partial_1)$. Since $-\Delta v=\nabla^{\perp}\omega$, this boundary condition follows by taking the normal trace of $(1.10)$. The problem (2.4) has a unique solution satisfying 

\begin{align*}
\sup_{x\in \Omega}d(x)|\nabla q(x)|\leq C||\omega||_{L^{\infty}(\partial\Omega)},    \tag{2.5}
\end{align*}\\
\cite{AG1}, \cite{AG2}, \cite{KLS} and by using the solution operator $\mathbb{K}: \omega\longmapsto \nabla q$, the associated pressure gradient is represented by $\nabla q=\mathbb{K}\omega$ for $v=R(\lambda)f$ and $f\in L^{\infty}_{\sigma}$.

\noindent
(ii) The operator $R(\lambda)$ is pseudo-resolvent on $L^{\infty}_{\sigma}$ with the trivial kernel, i.e. $\textrm{Ker}\ R(\lambda)=\{0\}$. Indeed, if $v=R(\lambda)f=0$, we have $\nabla q=\mathbb{K}\omega=0$ and $f=0$. Hence by the open mapping theorem, there exits a closed operator $A$ such that $R(\lambda)=(\lambda-A)^{-1}$. We call $A$ the Stokes operator on $L^{\infty}_{\sigma}$.
\end{rems}

\vspace{15pt}

We shall prove the representation formula (1.15) for solutions of (1.10) with the kernels (1.16) and (1.17).

\vspace{15pt}

\begin{lem}[Representation formula]
The solution $v=R(\lambda)f$ and $\nabla q=\mathbb{K}\omega$ for $\lambda\in \Sigma_{\theta+\pi/2}$ and $f\in L^{\infty}_{\sigma}$ is represented by 

\begin{align*}
v(x)=\int_{\Omega}E^{\lambda}(x-y)f(y)\dd y
+\int_{\partial\Omega}V^{\lambda}(x-y)TN(y)\dd H(y),\quad x\in \Omega,  \tag{2.6}
\end{align*}\\
for $T=\nabla v+{}^{t}\nabla v-qI$.
\end{lem}

\vspace{5pt}

\begin{proof}
We denote by $\overline{f}$ the zero extension of $f$ to $\mathbb{R}^{2}\backslash \overline{\Omega}$. Observe that $(\overline{v},\overline{q})$ is a weak solution of the problem 

\begin{equation*}
\begin{aligned}
\lambda \overline{v}-\Delta \overline{v}+\nabla \overline{q}=\overline{f}+\mu,\quad \D\ \overline{v}=0\quad \textrm{in}\ \mathbb{R}^{2},
\end{aligned}
\tag{2.7}
\end{equation*}\\
for a measure $\mu$ satisfying 

\begin{align*}
(\mu,\varphi)=\int_{\partial\Omega}TN(y)\cdot \varphi(y)\dd H(y),\quad \varphi\in C_{0}(\mathbb{R}^{2}),   \tag{2.8}
\end{align*}\\
where $C_{0}(\mathbb{R}^{2})$ denotes the space of all continuous functions in $\mathbb{R}^{2}$ vanishing at space infinity and $(\cdot,\cdot )$ denotes the pairing between $C_{0}(\mathbb{R}^{2})$ and its adjoint space. Indeed, multiplying $\varphi\in C^{\infty}_{c}(\mathbb{R}^{2})$ by (1.10) and integration by parts imply (2.7) in a weak sense. The formula (2.6) formally follows by multiplying $(\lambda-\Delta)^{-1}\mathbb{P}$ by (2.7). We set $v_1=(\lambda-\Delta)^{-1}\overline{f}$ and $v_2=\overline{v}-v_1$ to see that 

\begin{align*}
\lambda v_2-\Delta v_2+\nabla q=\mu,\quad \D\ v_2=0\quad \textrm{in}\ \mathbb{R}^{2}.
\end{align*}\\
By the mollifications $v_{2,\varepsilon}=v_2*\eta_{\varepsilon}$, $q_{\varepsilon}=q*\eta_{\varepsilon}$ and $\mu_{\varepsilon}=\mu*\eta_{\varepsilon}$ with the standard mollifier $\eta_{\varepsilon}$, $(v_{2,\varepsilon},q_{\varepsilon})$ satisfies the above problem for $\mu_{\varepsilon}\in L^{p}$ for $p\in [1,\infty]$. By multiplying $(\lambda-\Delta)^{-1}\mathbb{P}$ by the equation, we have

\begin{align*}
v_{2,\varepsilon}(x)=(\lambda-\Delta)^{-1}\mathbb{P}\mu_{\varepsilon}
=\int_{\mathbb{R}^{2}}V^{\lambda}(x-y)\mu_{\varepsilon}(y)\dd y
=\eta_{\varepsilon}*\left(\int_{\partial\Omega}V^{\lambda}(x-y)TN(y)\dd H(y)   \right).
\end{align*}\\
Sending $\varepsilon\to0$ yields (2.6). This completes the proof.
\end{proof}

\vspace{15pt}

The Stokes paradox follows a similar argument using the fundamental tensor of the Stokes equations. The following result is due to Chang and Finn \cite[Theorem 3]{ChangFinn}.

\vspace{15pt}

\begin{lem}[Stokes paradox]
Let $(v, \nabla q)\in W^{2,p}_{\textrm{loc}}(\overline{\Omega})\times L^{p}_{\textrm{loc}}(\overline{\Omega})$, $p\in (1,\infty)$, satisfy (1.13). Assume that 

\begin{align*}
v(x)=o(\log{|x|})\quad \textrm{as}\ |x|\to\infty.  \tag{2.9}
\end{align*}\\
Then, $v\equiv 0$ and $\nabla q\equiv 0$.
\end{lem}

\vspace{5pt}

\begin{proof}
We give a proof for completeness. Observe that the zero extension $(\overline{v}, \overline{q})$ is a solution of the problem 

\begin{equation*}
\begin{aligned}
-\Delta \overline{v}+\nabla \overline{q}=\mu,\quad \D\ \overline{v}=0\quad \textrm{in}\ \mathbb{R}^{2},
\end{aligned}
\tag{2.10}
\end{equation*}\\
for a measure $\mu$ defined by (2.8). By the fundamental tensor of the Stokes equations $V=(V_{ij})$ and $Q=(Q_{j})$ \cite[p.239]{Gal}, 

\begin{align*}
V_{ij}(x)=\frac{1}{4\pi}\left(-\delta_{ij}\log{|x|}+\frac{x_ix_j}{|x|^{2}}\right),\quad Q_{j}(x)=\frac{1}{2\pi}\frac{x_j}{|x|^{2}},
\end{align*}\\
we set $(\tilde{v}, \tilde{q})$ by 

\begin{align*}
\tilde{v}(x)=\int_{\partial\Omega}V(x-y)TN(y)\dd H(y),\quad  \tilde{q}(x)=\int_{\partial\Omega} Q(x-y)\cdot TN(y)\dd H(y).   
\end{align*}\\
The functions $\tilde{v}$ and $\tilde{q}$ are locally integrable in $\mathbb{R}^{2}$ and $\tilde{v}=O(\log{|x|})$, $\nabla \tilde{v}=O(|x|^{-1})$ as $|x|\to\infty$. Observe that $u=\overline{v}-\tilde{v}$ and $p=\overline{q}-\tilde{q}$ is a weak solution of 

\begin{align*}
-\Delta u+\nabla p=0,\quad \D\ u=0\quad \textrm{in}\ \mathbb{R}^{2}.  \tag{2.11}
\end{align*}\\
Since $u$ and $p$ are locally integrable in $\mathbb{R}^{2}$, by mollification we may assume that they are smooth in $\mathbb{R}^{2}$. Since $\omega=\partial_1u^{2}-\partial_2 u^{1}$ is bounded in $\mathbb{R}^{2}$ and satisfies $-\Delta \omega=0$ in $\mathbb{R}^{2}$, $\omega$ is constant by the Liouville theorem. By $-\Delta u=0$ in $\mathbb{R}^{2}$ and $u=O(\log{|x|})$ as $|x|\to\infty$, $u$ and $p$ are constants. Hence by shifting the pressure up to constant 

\begin{align*}
v(x)=v_{\infty}+\int_{\partial\Omega}V(x-y)TN(y)\dd H(y),\quad q(x)=\int_{\partial\Omega} Q(x-y)\cdot TN(y)\dd H(y),   \tag{2.12}
\end{align*}\\
for some constant $v_{\infty}$. This implies 

\begin{align*}
v(x)&=v_{\infty}+V(x)\int_{\partial\Omega}TN(y)\dd H(y)+O(|x|^{-1}),\\
q(x)&=Q(x)\cdot \int_{\partial\Omega}TN(y)\dd H(y)+O(|x|^{-2})\quad \textrm{as}\ |x|\to\infty.
\end{align*}\\
Since $v=o(\log{|x|})$ as $|x|\to\infty$, by dividing $v$ by $\log{|x|}$ and sending $|x|\to\infty$, 

\begin{align*}
\int_{\partial\Omega}TN(y)\dd H(y)=0.
\end{align*}\\
Hence $v-v_{\infty}=O(|x|^{-1})$ and $\nabla v,\ q=O(|x|^{-2})$ as $|x|\to\infty$. By multiplying $v-v_{\infty}$ by (1.13) and integration by parts in $\Omega\cap B(0,R)$,

\begin{align*}
\int_{\Omega\cap B(0,R)}|\nabla v|^{2}\dd x=\int_{\partial B(0,R)}(TN)\cdot (v-v_{\infty})\dd H(x)\to 0\quad \textrm{as}\ R\to\infty.
\end{align*}\\
By $v=0$ on $\partial\Omega$, $v\equiv 0$ and $\nabla q\equiv 0$ follow. This completes the proof.
\end{proof}

\vspace{15pt}

\begin{rem}
For $n\geq 3$, the fundamental tensor of the Stokes equations (2.11) is $V=(V_{ij})$, $Q=(Q_j)$ for 

\begin{align*}
V_{ij}(x)=\frac{1}{2n(n-2)\alpha(n)}\left(\frac{\delta_{ij}}{|x|^{n-2}}+(n-2)\frac{x_ix_j}{|x|^{n}}\right),\quad Q_j(x)=\frac{1}{n\alpha(n)}\frac{x_j}{|x|^{n}}.
\end{align*}\\ 
In the same way as the proof of Lemma 2.5, we see that any bounded solutions $v$ of (1.13) is of the form (2.12) for some constant $v_{\infty}$.
\end{rem}

\vspace{15pt}

\section{The resolvent estimate}

\vspace{15pt}

We prove the estimate (1.11). By the approximation for $f\in L^{\infty}_{\sigma}$ (Lemma 2.2), it suffices to show (1.11) for $f\in C_{c,\sigma}^{\infty}$.

\vspace{15pt}

\begin{prop}
\begin{align*}
\sup_{\lambda\in \Sigma_{\theta+\pi/2}}|\lambda|\ ||R(\lambda)f||_{L^{\infty}}<\infty,\quad f\in C_{c,\sigma}^{\infty}.   \tag{3.1}
\end{align*}
\end{prop}

\vspace{5pt}

\begin{proof}
Since $C_{c,\sigma}^{\infty}\subset L^{p}_{\sigma}$ for $p\in (1,\infty)$, $R(\lambda)f=(\lambda-A)^{-1}f$ for $A=\mathbb{P}\Delta$ and the Helmholtz projection operator $\mathbb{P}$. The domain $D(A)=W^{2,p}\cap W^{1,p}_{0}\cap L^{p}_{\sigma}$ is equipped with the graph-norm and $D(A)\subset W^{2,p}$ with continuous injection \cite{G81}. Here, $W^{2,p}$ denotes the Sobolev space and $W^{1,p}_{0}$ denotes the space of all trace zero functions in $W^{1,p}$. By the $L^{p}$-resolvent estimate $|\lambda|\ ||R(\lambda)f||_{L^{p}}\leq C||f||_{L^{p}}$ \cite{BV} and the Sobolev embedding for $p\in (2,\infty)$,

\begin{align*}
||R(\lambda)f||_{L^{\infty}}
\leq C||R(\lambda)f||_{W^{2,p}}
\leq C'\left(||R(\lambda)f||_{L^{p}}+||AR(\lambda)f||_{L^{p}} \right) 
\leq C''\left(\frac{1}{|\lambda|}+1\right)||f||_{L^{p}}.
\end{align*}\\
Hence $|\lambda |\ ||R(\lambda)f||_{L^{\infty}}$ is bounded for $|\lambda|\leq 1$. Since $|\lambda|\ ||R(\lambda)f||_{L^{\infty}}\leq C||f||_{L^{\infty}}$ for $|\lambda|\geq 1$ by (2,2), (3.1) follows.
\end{proof}

\vspace{15pt}

\begin{lem}
There exists a constant $C$ such that (1.11) holds for $f\in C_{c,\sigma}^{\infty}$ and $\lambda\in \Sigma_{\theta+\pi/2}$.
\end{lem}

\vspace{5pt}

\begin{proof}
We argue by contradiction. Suppose that (1.11) were false. Then, for $m\geq 1$ there exists $\tilde{f}_m\in  C^{\infty}_{c,\sigma}$ such that 

\begin{align*}
M_m=\sup_{\lambda\in \Sigma_{\theta+\pi/2}}  |\lambda|\ ||R(\lambda)\tilde{f}_m||_{L^{\infty}}(\lambda)> m||\tilde{f}_m||_{L^{\infty}}.
\end{align*}\\
By setting $f_m=\tilde{f}_m/M_m$, 

\begin{align*}
\sup_{\lambda\in \Sigma_{\theta+\pi/2}}  |\lambda|\ ||R(\lambda)f_m||_{L^{\infty}}(\lambda)=1,\quad  ||f_m||_{L^{\infty}}< \frac{1}{m}.
\end{align*}\\
We set $v_m=R(\lambda)f_m$ and take a point $\lambda_m\in \Sigma_{\theta+\pi/2}$ such that 

\begin{align*}
|\lambda_m|\ ||v_m||_{L^{\infty}}\geq \frac{1}{2}.
\end{align*}\\
We may assume that $\lambda_m\to0$ by (2.2). Observe that $u_m=\lambda_m v_m$ satisfies 

\begin{align*}
\lambda_m u_m-\Delta u_m+\nabla p_m=\lambda_m f_m,\quad \D\ u_m=0\quad \textrm{in}\ \Omega,\\
u_m=0\quad \textrm{on}\ \partial\Omega,
\end{align*}\\
with some associated pressure $p_m$. We take a point $x_m\in \Omega$ such that 

\begin{align*}
|u_m(x_m)|\geq \frac{1}{4}.
\end{align*}\\
We normalize the pressure $p_m$ so that $\int_{\partial\Omega}p_m\dd H(y)=0$. Since $u_m-\Delta u_m+\nabla p_m=\lambda_m(f_m-u_m)+u_m$, applying the resolvent estimates (2.1) for $p>2$ implies

\begin{align*}
||u_m||_{W^{1,\infty}}+\sup_{x\in \Omega}\left\{||\nabla^{2} u_m||_{L^{p}(\Omega\cap B(x,1) ) }+||\nabla p_m||_{L^{p}(\Omega\cap B(x,1) ) }\right\}&\leq C(||\lambda_m(f_m-u_m)||_{L^{\infty}}+||u_m||_{L^{\infty}}  ) \\
&\leq C',\quad  \textrm{for all}\ m\geq 1.
\end{align*}\\
Hence $\{u_m\}$ is equi-continuous in $\overline{\Omega}$. By choosing a subsequence (still denoted by $\{u_m\}$), $u_m$ converges to a limit $u$ locally uniformly in $\overline{\Omega}$ together with $\nabla u_m$ and $p_m$. Then the limit $u$ is a bounded solutions of 

\begin{align*}
-\Delta u+\nabla p=0,\quad \D\ u=0\quad \textrm{in}\ \Omega,\\
u=0\quad \textrm{on}\ \partial\Omega,
\end{align*}\\
with the associated pressure $p$. Applying Lemma 2.5 implies that $u\equiv 0$ and $\nabla p\equiv 0$. Since $\int_{\partial\Omega}p\dd H(y)=0$, $p\equiv 0$. Hence we have 

\begin{align*}
u_m\to 0\quad \textrm{locally uniformly in}\ \overline{\Omega},  \tag{3.2}
\end{align*}\\
together with $\nabla u_m$ and $p_m$. In particular, $T_m=\nabla u_m+{}^{t}\nabla u_m-p_mI\to 0$ uniformly on $\partial\Omega$ as $m\to\infty$.

Suppose that $\lim\sup_{m\to\infty}|x_m|<\infty$. By choosing a subsequence, we may assume that $\{x_m\}$ converges to some point in $\overline{\Omega}$. This implies that $1/4\leq |u_m(x_m)|\to0$, a contradiction. We may assume that $\lim\sup_{m\to\infty}|x_m|=\infty$. By choosing a subsequence, we may assume that $\lim_{m\to\infty}|x_m|=\infty$. We consider two cases depending on whether $|\lambda_m|^{1/2}|x_m|$ vanishes or not.\\

\noindent 
\noindent 
\textit{Case 1}. $\lim\inf_{m\to\infty}|\lambda_m|^{1/2} |x_m|>0$.\\

We may assume that $|\lambda_m|^{1/2} |x_m|\geq  d$ for some constant $d>0$ by choosing a subsequence. By the representation formula (2.6), 

\begin{align*}
u_m(x)=(\lambda_m-\Delta)^{-1}\lambda_m f_m+\int_{\partial\Omega}V^{\lambda_m}(x-y)T_mN(y)\dd H(y).
\end{align*} \\
By $|\lambda_m|^{1/2}|x_m|\geq d$ and (1.18),

\begin{align*}
\sup_{y\in \partial\Omega}|V^{\lambda_m}(x_m-y)|\leq C,\quad \textrm{for all}\ m\geq 1.
\end{align*}\\
By the $L^{\infty}$-estimate $|\lambda_m|\ ||(\lambda_m-\Delta)^{-1}f_m||_{L^{\infty}}\leq C||f_m||_{L^{\infty}}$, 

\begin{align*}
\frac{1}{4}\leq |u_m(x_m)|\leq C\left(\frac{1}{m}+\int_{\partial\Omega}|T_mN(y)|\dd H(y)\right)\to 0\quad \textrm{as}\  m\to\infty.
\end{align*}\\
Thus Case 1 does not occur. \\

\noindent 
\noindent 
\textit{Case 2}. $\lim\inf_{m\to\infty}|\lambda_m|^{1/2} |x_m|=0$.\\

We may assume that $\lim_{m\to\infty}|\lambda_m|^{1/2} |x_m|= 0$. By the representation formula (2.6) and (1.19),

\begin{equation*}
\begin{aligned}
u_m(x)&=(\lambda_m-\Delta)^{-1}\lambda_m f_m
-\frac{1}{4\pi}\log{\sqrt{\lambda_m}}\int_{\partial\Omega} T_m N(y)\dd H(y) \\
&-\frac{1}{4\pi}\int_{\partial\Omega}\log{|x-y|}T_m N(y)\dd H(y)
+\int_{\partial\Omega}\tilde{V}^{\lambda_m}(x-y)T_mN(y)\dd H(y).
\end{aligned}
\tag{3.3}
\end{equation*}\\
By (3.2), sending $m\to\infty$ for fixed $x\in \Omega$ implies 

\begin{align*}
0=\lim_{m\to\infty}\log|\lambda_m|^{1/2} \left|\int_{\partial\Omega}T_mN(y)\dd H(y)\right|.   \tag{3.4}
\end{align*}\\
We substitute $x=x_m$ into (3.3). By (1.18),

\begin{align*}
\sup_{y\in \partial\Omega}|\tilde{V}^{\lambda_m}(x_m-y)|\leq C,\quad \textrm{for all}\ m\geq 1.
\end{align*}\\
Since $|\lambda_m|^{1/2}|x_m|\leq 1$ for sufficiently large $m$, $\log|x_m|\leq -\log |\lambda_m|^{1/2}$ and 

\begin{align*}
\left|\int_{\partial\Omega}\log|x_m-y| TN(y)\dd H(y)\right|
&\leq \int_{\partial\Omega}\left|\log\left|\frac{x_m}{|x_m|}-\frac{y}{|x_m|}\right|\right| |T_mN(y)|\dd H(y)
+\log|x_m| \left|\int_{\partial\Omega} T_mN(y)\dd H(y)\right| \\
&\leq C\int_{\partial\Omega}|T_mN(y)|\dd H(y)
-\log|\lambda_m|^{1/2}  \left|\int_{\partial\Omega} T_mN(y)\dd H(y)\right|.
\end{align*}\\
By (3.4) and the dominated convergence theorem,

\begin{align*}
\frac{1}{4}\leq |u_m(x_m)|
&\leq \frac{C}{m}-\frac{1}{2\pi}\log|\lambda_m|^{1/2}  \left|\int_{\partial\Omega} T_mN(y)\dd H(y) \right| 
+C\int_{\partial\Omega}\left|T_mN(y)\right|\dd H(y)\\
&\to 0\quad \textrm{as}\ m\to\infty.
\end{align*}\\
We obtained a contradiction. Thus Case 2 does not occur.

We conclude that both Case 1 and Case 2 do not occur. The proof is now complete.
\end{proof}

\vspace{15pt}

\begin{proof}[Proof of Theorem 1.3]
For $f\in L^{\infty}_{\sigma}$, we take a sequence $\{f_m\}\subset C_{c,\sigma}^{\infty}$ satisfying (2.3) by Lemma 2.2 (i). Since $|\lambda|\ ||R(\lambda)f_m||_{L^{\infty}}\leq C||f||_{L^{\infty}}$ for all $m\geq 1$ and $R(\lambda)f_m$ converges to $R(\lambda)f$ locally uniformly in $\overline{\Omega}$ by Lemma 2.2 (ii), the limit satisfies the desired estimate. Hence the assertion (i) holds. The assertion (ii) follows from the Dunford integral of the resolvent by using (1.11).
\end{proof}

\vspace{15pt}

\begin{rems}
(i) Besides the estimate (1.12), we obtain estimates for spatial derivatives, 

\begin{align*}
||\nabla S(t)v_0||_{L^{\infty}}+||\nabla^{2} S(t)v_0||_{L^{\infty}}\leq C||v_0||_{L^{\infty}},\quad t\geq 1, \quad v_0\in L^{\infty}_{\sigma}.  \tag{3.5}
\end{align*}\\
This follows from (1.12) and the finite time estimate $t^{1/2}||\nabla S(t)v_0||_{L^{\infty}}+t||\nabla^{2} S(t)v_0||_{L^{\infty}}\leq C||v_0||_{L^{\infty}}$ for $0<t\leq T$ \cite{AG2}.

\noindent
(ii) For $n\geq 2$ and $v_0\in L^{\infty}_{\sigma}$, Lemma 2.5 implies that 

\begin{align*}
S(t)v_0\to 0\quad \textrm{locally uniformly in}\ \overline{\Omega}\ \textrm{as}\ t\to\infty.  \tag{3.6}
\end{align*}\\
In fact, suppose that (3.6) were false. Then, there exists a sequence $\{t_m\}$ such that $t_m\to\infty$ and (3.6) does not hold. By (1.12), (3.5) and choosing a subsequence (still denoted by $\{t_m\}$) $v_m(t)=S(t+t_m)v_0$ converges to a limit $v$ locally uniformly in $\overline{\Omega}\times [0,\infty)$. Since the limit $v$ is bounded and independent of $t$, $v\equiv 0$ by Lemma 2.5 and $S(t_m)v_0\to 0$ locally uniformly in $\overline{\Omega}$. This is a contradiction. 

\noindent
(iii) For $n\geq 3$ and $v_0\in L^{\infty}_{\sigma}$, 

\begin{align*}
S(t)v_0\to v\quad \textrm{locally uniformly in}\ \overline{\Omega}\ \textrm{as}\ t\to\infty,  \tag{3.7}
\end{align*}\\
for some solution $v$ of the stationary Stokes equations (1.13). Since any bounded solutions of (1.13) for $n\geq 3$ must be asymptotically constant as $|x|\to\infty$ by Remark 2.6, $S(t)v_0$ is asymptotically constant as $t\to\infty$ and $|x|\to\infty$ for any bounded initial data $v_0\in L^{\infty}_{\sigma}$.
\end{rems}

\vspace{15pt}

\section*{Acknowledgements}
The author is partially supported by JSPS through the Grant-in-aid for Young Scientist (B) 17K14217, Scientific Research (B) 17H02853 and Osaka City University Advanced Mathematical Institute (MEXT Joint Usage / Research Center on Mathematics and Theoretical Physics).

\vspace{15pt}

\section*{Conflict of Interest}
The author declares that he has no conflict of interest.

\vspace{15pt}

\bibliographystyle{abbrv}
\bibliography{ref}

\end{document}